\documentclass{article}
\usepackage{graphicx}
 \usepackage{mathptmx}
\usepackage{amsmath, amstext,amssymb,amsfonts}
\usepackage[english]{babel}
\usepackage{delarray}
\usepackage{mathptmx}
\usepackage{amsmath, amstext,amssymb,amsfonts}
\usepackage[english]{babel}
\usepackage{delarray}

\DeclareFontFamily{U}{mathx}{\hyphenchar\font45}
\DeclareFontShape{U}{mathx}{m}{n}{
      <5> <6> <7> <8> <9> <10>
      <10.95> <12> <14.4> <17.28> <20.74> <24.88>
      mathx10
      }{}
\DeclareSymbolFont{mathx}{U}{mathx}{m}{n}
\DeclareFontSubstitution{U}{mathx}{m}{n}
\DeclareMathAccent{\widecheck}{0}{mathx}{"71}
\DeclareMathAccent{\wideparen}{0}{mathx}{"75}

\numberwithin{equation}{section}

\newtheorem{thm}{Theorem}[section]
 
 \newtheorem{lem}[thm]{Lemma}
 \newtheorem{prop}[thm]{Proposition}

\newcommand{\Co}{\mathbb C} \newcommand{\R}{\mathbb R}   \newcommand{\Z}{\mathbb Z}
\newcommand{\Zn}{\mathbb{Z}^n}     \newcommand{\Rn}{\mathbb{R}^{n}}  \newcommand{\Con}{\mathbb{C}^{n}}
\newcommand{\qun}{Q^{n}_{\sigma}} \newcommand{\qup}{Q^{n}_{\pi}}     \newcommand{\bnp}{B^p_{\qup}}     \newcommand{\bnn}{B_{\qup}} \newcommand{\bp}{B_{\qun}^{p}}\newcommand{\be}{B_{\sigma}^{p}}  \newcommand{\M}{M({\mathbb R}^n)}
\newcommand{\bi}{B^{\infty}_{\qun}}            \newcommand{\bip}{B^{\infty}_{Q^{n}_{\pi}}}
 \newcommand{\sic}{{\rm{sic}}_{n}}  \newcommand{\sinc}{{{\rm{sinc}}_{n}}\,} \newcommand{\sivc}{{{\rm{sinc}}_{1}}\,}

\date{}
\title{On a sampling expansion with  partial derivatives for functions of several variables}
\author{{\large Saulius  Norvidas} }
\date{{\footnotesize Institute of Data Science and Digital Technologies, Vilnius University, Akademijos str. 4, Vilnius LT-04812, Lithuania\\
 ({\rm{e-mail: norvidas{@}gmail.com}})}}
\begin{document}

\maketitle
 {{ {\bf Abstract}}}
 Let $B^p_{\sigma}$, $1\le p<\infty$, $\sigma>0$,  denote the space of all $f\in L^p(\R)$ such that the Fourier transform of $f$ (in the sense of distributions) vanishes outside $[-\sigma,\sigma]$. The classical sampling theorem states  that each  $f\in B^p_{\sigma}$  may be reconstructed exactly from  its sample values  at equispaced sampling points $\{\pi m/\sigma\}_{m\in\Z} $ spaced by $\pi /\sigma$.  Reconstruction is also possible  from   sample values  at  sampling points $\{\pi \theta m/\sigma\}_m $ with certain $1< \theta\le 2$ if we  know    $f(\theta\pi m/\sigma) $ and  $f'(\theta\pi m/\sigma)$, $m\in\Z$. In this paper we present  sampling series for functions of several variables. These series involves  samples of functions  and    their partial  derivatives.



\section{Introduction}\label{s:1}
{\large{

\ \ We start with some notation and definitions. Let $\Zn$, $\Rn$ and $\Con$ be the $n$-dimensional integer lattice,  the real  Euclidean space and the complex Euclidean space, respectively. For any $\tau\in\Co$ and each $a,b\in\Con$  we write
\[
\tau a=(\tau a_1,\dots, \tau a_n),\quad ab=(a_1b_1,\dots,a_nb_n).
\]
If, in addition, $b_j\neq 0$, $j=1,\dots,n$, then $a/b$ denotes the  vector of fractions $\bigl(a_1/b_1,\dots, a_n/b_n\bigr)$. For  $\sigma\in\Rn$  such that $\sigma_j>0$, $j=1,\dots,n$,  let us denote by $\sigma\Zn$ the lattice $\oplus_{j=1}^n \sigma_j\,\Z$. Also if $A, B\subset \Con$, then $A+B=\{a+b: a\in A, b\in B\}$.

\ \ For $f\in L^1(\Rn)$, we define the Fourier of $f$ transform    by
\[
\widehat{f}(t)= \int_{\Rn} e^{-i\langle x,t\rangle}f(x)\,dx,
\]
 $ x\in\Rn$, where $\langle  x, t\rangle =\sum_{k=1}^{n}  x_kt_k$ is the  scalar product on $\Rn$.  If $f\not\in L^{1}(\Rn)$,  then  we understand  $\widehat{f}$ in a distributional sense of   tempered distributions $S'(\Rn)$. Given a closed subset $\Omega\subset\Rn$,  a function ƒ$f:\Rn\to \Co$     is called bandlimited to ƒ$\Omega$  if ƒ$\widehat{f}$  vanishes outside $\Omega$.

\ \  For $1\le p\le\infty$ and $\sigma\in\Rn$ such that  ƒ$\sigma_j>0$, $j=1,\dots, n$, let
 \[
 \qun=\{x\in\Rn: |x_j|\le \sigma_j, j=1,\dots,n\}\quad{\text{\rm{and}}}\quad \bp=\{f\in L^p(\Rn):\  supp \ \widehat{f}\subset \qun\}.
\]
The space $\bp$  is a Banach space in  $L^p(\Rn)$-norm. In the case of functions of one variable, let us write $\be$  instead of $\bp$.

  \ \ By the Paley-Wiener-Schwartz theorem (see \cite[p. 181]{7}), each   $f\in \bp$, $1\le p\le\infty$,  is infinitely differentiable on $\Rn$. Moreover, $f$ has   an extension onto  $\Con$  to an entire function. Note that we shall identify any $f\in \bp$ with a $L^p(\Rn$)-function $f(x)$ defined on $\Rn$ and in other cases consider the same $f$ as entire function $f(z)$ defined on the whole $\Con$. In the sequel, we shall frequently use the following  functions defined on $\Con$
\[
\sinc(z)=\prod_{j=1}^{n}\frac{\sin\pi z_j}{\pi z_j}\quad{\text{\rm{and}}}\quad \sic(z)=\prod_{j=1}^{n}\sin\frac{ z_j}2.
\]
Of course, $\sic\in\bi$ with $\sigma=(1/2,\dots,1/2)$ and $\sinc\in \bnp$ for each $1< p\le\infty$.

\ \ The classical Whittaker-Shanon-Kotelnikov  theorem states that for any $f\in B^p_{\sigma}$, $1\le p<\infty$,  the following sampling expansion holds (see, e.g., \cite[p. 51]{6})
\begin{equation}
 f(z)= \sum_{m\in\Z}f\Bigl(\frac{\pi\, m}{\sigma}\Bigr){{{\rm{sinc}}_{1}}\,}\Bigl(\frac{\sigma}{\pi}\Bigl(z-\frac{\pi\,m}{\sigma}\Bigr)\Bigr)=
 \sum_{u\in(\pi/\sigma)\Z} f(u){{\rm{sinc}}}_1\Bigl(\frac{\sigma}{\pi}(z-u)\Bigr).
 \end{equation}
In particular, this   series  converges absolutely for $z\in\Co$ and uniformly on $\R$ and also on compact subsets of $\Co$. Hence,   (1.1) shows that each $f\in B^p_{\sigma}$, $1\le p<\infty$, may be reconstructed exactly from  sample values $f(\pi m/\sigma) $ at equispaced sampling points $\{\pi m/\sigma\}_{m\in\Z} $ spaced by $\pi /\sigma$.  Note, that this spaced value $\pi /\sigma$ is exact, i.e., for any $\theta>1$ there exist two $f_{j}\in B^p_{\sigma}$, $j=1,2$ such that  $f_1\not\equiv f_2$, but  $f_1(\theta\pi m/\sigma)= f_2(\theta\pi m/\sigma)$ for all $m\in\Z$. On the other hand, if we  know  sample   values  $f(\theta\pi m/\sigma) $ with certain $1< \theta\le 2$,  then the reconstruction of $f$ is   possible in the case if we also use its derivative values $f'(\theta\pi m/\sigma) $. In particular, if $f\in\be$ with $1\le p<\infty$, then  (see e.g.,  \cite[p. 145]{8})
\begin{equation}
f(z)=\sum_{u\in(2\pi/\sigma)\Z}\Bigl(f(u)+f'(u)(z-u)\Bigr){{\rm{sinc}}}^2_1\Bigl(\frac{\sigma}{2\pi}(z-u)\Bigr).
 \end{equation}

\ \ The following $n$-dimensional sampling theorem is  a standard  extension of  (1.1) in $\bp$, $1\le p<\infty$  (see, e.g.,  \cite[p. 172]{5})
 \begin{equation}
 f(z)=\sum_{u\in(\pi/\sigma)\Zn} f(u)\sinc\Bigl(\frac{\sigma}{\pi}(z-u)\Bigr).
   \end{equation}
\ \   The aim of this paper is to prove a multidimensional version of (1.2). Note that in \cite{2}  (see also  a tutorial review [9, p. 40]) was given the following expression
\begin{gather}
f(z)=\sum_{u\in\Z^2}\biggl[ f\Bigl(2\pi \frac{u}{\sigma}\Bigr)+\Bigl(z_1-\frac{2\pi u_1}{\sigma_1}\Bigr)\cdot\frac{\partial f}{\partial z_1}\Bigl(2\pi \frac{u}{\sigma}\Bigr)+\Bigl(z_2-\frac{2\pi u_2}{\sigma_2}\Bigr)\cdot\frac{\partial f}{\partial z_2}\Bigl(2\pi \frac{u}{\sigma}\Bigr)\biggr]\nonumber\\
{{{\rm{sinc}}_{2}}\,}^2\Bigl(\frac1{2\pi}\sigma\Bigl(z-u)\Bigr)
\end{gather}
for  $f\in B^{p}_{Q^2_{\sigma}}$ with $1\le p<\infty$.  We say that such a sampling theorem  fails in general. Indeed, let $\chi$ be any function in the Schwartz space $S(\R^2)$ such that $\chi\not\equiv 0$ and  ${\text{supp}}\, \chi\subset Q^2_{\sigma/2}$. Then the function
\begin{equation}
f(z)= \widehat{\chi}(z) {{{\rm{sic}}_{2}}\,}(\sigma z)
 \end{equation}
is in $B^{p}_{Q^2_{\sigma}}$ for each $ 1\le p\le \infty$. Moreover,
\[
f\Bigl(2\pi\frac{u}{\sigma}\Bigr)=\frac{\partial f}{\partial  z_1}\Bigl(2\pi\frac{u}{\sigma}\Bigr)=\frac{\partial f}{\partial  z_2}\Bigl(2\pi\frac{u}{\sigma}\Bigr)=0
\]
for all $u\in\Z^2$. Hence, in this case the series (1.4) generates the zero function, but not the function (1.5). Even more, (1.5) shows that the same still true if we added to (1.4) an arbitrary number of the following sample values
\[
\frac{\partial^{m} f}{\partial  x_j^m}\Bigl(2\pi\frac{u}{\sigma}\Bigr)
\]
with $u\in\Z^2$, $j=1,2$ and $m=2,3,\dots$. Therefore, any  multidimensional version of (1.2) must necessarily  contains also  mixed   partial derivatives of $f$.

\ \ Now we shall provide some more notation and formulate our main theorem. For $k\in\Zn$ with $k_j\ge 0$, $j=1,\dots,n$, and $z\in\Con$, here and subsequently, we denote  the operator
\[
\frac{\partial^{|k|}}{\partial^{k_1}z_1\dots\partial^{k_n}z_n}, \quad |k|=k_1+\dots+k_n
\]
by $\partial^k_{z}$ for short. Note that if $k_j=0$ for all  $j=1,\dots,n$, then $\partial^k_{z} f(z)$ is simply $f(z)$. Set
\[
E^n=\bigl\{(t_1,\dots,t_n): \ t _j\in\{0;1\}, j=1,\dots, n\bigr\}.
\]
For fixed $k\in E^n$, $v\in\Zn$ and $f\in\bp$, let us define the following polynomial in $\lambda\in\Con$ by
\begin{equation}
P_{f,k,v}(\lambda)=\Bigl(\partial ^k_z f\Bigr)(v)\cdot\prod_{j=1}^n \lambda^{k_j}_j.
 \end{equation}
\begin{thm}.
 Let $f\in\bp$ with $1\le p<\infty$. Then
\begin{equation}
f(z)=\sum_{u\in(2\pi/\sigma)\Zn}\biggl(\sum_{k\in E^n}P_{f,k,u}(z-u)\biggr)\sinc^2\Bigl(\frac1{2\pi}\sigma(z-u)\Bigr).
\end{equation}
The series (1.7) converges absolutely and uniformly on $\Rn$ and also on any compact subset of $\Con$.
\end{thm}
\  For fixed   $u\in(2\pi/\sigma)\Zn$, the representation (1.7) contains $2^n$ values of partial derivatives $\partial ^k_z f(u)$ when $k$ obtain all possible values from $  E^n$. Note that  this   representation   is exact in some sense. More precisely, if we eliminate  in (1.7) an arbitrary  polynomial $P_{f,k,u}$ with certain $k=\widetilde{k}\in E^n$, then such a formula will be false. Indeed, let
\[
\widetilde{f}(z)= \widehat{\chi}(z)\frac{\partial^{n-|\widetilde{k}|}  {{{\rm{sic}}_{n}}\,}(\sigma z)}
{\partial^{1-\widetilde{k}_1}z_1\dots\partial^{1-\widetilde{k}_n}z_n }
\]
where $\widehat{\chi}(z)$ is the same as that  of (1.5). Then an easy computation shows that
\[
\partial^k_z\widetilde{f}\Bigl(\frac{2\pi u}{\sigma}\Bigr)=0
\]
for all $u\in\Zn$ and each $k\in E^n$ such that $k\neq \widetilde{k}$. Therefore, in  such a  case  $\widetilde{f}$ will generate by (1.7) the zero function, but not $\widetilde{f}$.

\section{ Preliminaries and proofs}

\ \ In the sequel, we consider only the case $\qun=\qup=\{x\in\Rn: \max_{1\le j\le n}|x_j|\le \pi\}$. This involves no loss of generality, since the operator
\[
T_{\sigma}f(z)= \theta\cdot f\Bigl(\frac{\pi z}{\sigma}\Bigr)\quad{\text{with}}\quad \theta=\biggl(\frac{\prod_{j=1}^n\sigma_j}{\pi^n}\biggr)^{1/p},
\]
$z\in\Con$, is an isometric isomorphism  between $\bp$ and $\bnp$ for all $1\le p\le\infty$.

 \ \ For technical reasons, let us define  the following Banach space
 \[
 \bnn=\{f\in C_0(\Rn):\ supp \ \widehat{f}\subset \qup\},
 \]
  where  $C_0(\Rn)$ is the usual space of continuous functions on $\Rn$ that vanish at infinity.    If $1\le p<\infty$, then any $f\in \bnp$  satisfies $\lim_{|x|\to\infty} f(x)=0$ (see \cite[p. 118]{10}).  Next, if $1\le p< q<\infty$, then
 \begin{equation}
 \bnp\subset B_{\qup}^q\subset\bnn\subset B^{\infty}_{\qup}.
 \end{equation}

\ \ Given $m\in\{1,2,\dots,n\}$, set
\[
H^n_m=\{z=(z_1,\dots,z_n)\in\Con: z_m\in 2\Z\}.
\]
It is clear that
\begin{equation}
2k+H^n_m=H^n_m
 \end{equation}
 for each $k\in\Zn$. Let
\[
H^n=\bigcup_{m=1}^n H^n_{m}.
\]
 If $\varrho$ is a permutation of the set $\{1,2,\dots,n\}$, then
\begin{equation}
w\in H^n\quad {\rm{if \ and \ only \ if }}\quad (w_{\varrho(1)},w_{\varrho_(2)},\dots, w_{\varrho(n)})\in H^n.
 \end{equation}

\ \ Assume that $f:\Con\to\Co$ is analytic in a neighbourhood $U_a$ of $a\in\Con$ and $f(a)=0$. Let $f(z)=\sum_{m=0}^{\infty} P_m(z-a)$ be the expansion of $f$ into homogeneous polynomials in $(z-a)$-powers. Recall that  the minimal value of $m$ such that $P_m\not\equiv 0$ on $U_a$ is called the order of the zero $a$ for $f$. We denote this order by ${\rm{ord}}_a(f)$. Note that if $f(a)\neq 0$, then say that ${\rm{ord}}_a(f)=0$.

\begin{lem}
Let $f\in \bnn$. If $f(z)=0$ for all $z\in H^n$, then there is an entire function $ g:\Con\to\Co$ such that
\[
f(z)=\sic(\pi z)g(z),
\]
 $z\in\Con$.
\end{lem}
\ \ {\bf Proof}.\ It is easy to see that  the zeros set (all complex zeros) of $\sic(\pi z)$ coincides with  $H^n$. Therefore, the statement of our lema  follows immediately from  application to $f$ of the following fact (see  \cite[p. 12]{1}): if $F$ and $H$ are entire functions on $\Con$ such that ${\rm{ord}}_z(H)\le {\rm{ord}}_z(F)$ for all $z\in\Con$,  then there is an entire function $G$ such that $F\equiv GH$. The proof is complete.

\ \   Let  $M(\Rn)$ denote  the  Banach algebra of bounded regular   Borel measures  on  $\Rn$  with the total variation norm $\|\mu\|=\|\mu\|_{\M}$ and convolution as multiplication. The Fourier-Stieltjes transform of $\mu\in \M$   is given by
\[
\hat{\mu}(x)= \int_{\Rn} e^{-i\langle x,t\rangle}\, d\mu(t),\quad x\in\Rn.
\]
\ \ We need certain   facts about  differential and convolution operators  on $\bnn$. Bernstein's inequality   (see  \cite[p. 116]{10}) states that each partial derivative operator acts on $\bnp$, $1\le p\le \infty$, as  bounded operator. We do not find the  proof of the similar fact in the case of $\bnn$. For this reason, the proof of the following lemma is added here for completeness.
\begin{lem}
Let $f\in \bnn$.  Then $\partial^k_z f\in\bnn$ for all $k\in\Zn$ such that $k_1,\dots,k_n\ge 0$.
\end{lem}
\ \ {\bf Proof}.\ Let $\mu\in\M$. Then
\begin{equation}
T_{\mu}f(x)=\int_{\Rn}f(x-y)\,d\mu(y)=f\ast\mu(x),
 \end{equation}
$x\in\Rn$, is well-defined linear bounded operator operator on $\bip$    (see   \cite[p. 646]{3}). Next, if $\nu\in\M$  is such that $\widehat{\mu}=\widehat{\nu}$ on $\qup$, then $T_{\mu}=T_{\nu}$ on $\bip$ (\cite[p. 90]{4}). From (2.4) it follows that $\|T_{\mu}\|_{\bip}\le \|\mu\|_{\M}$. Moreover  (\cite[p. 646]{3}),
\begin{equation}
\|T_{\mu}\|_{\bip}=\inf\{\|\nu\|: \nu\in\M, \widehat{\nu}=\widehat{\mu} \ {\rm{on}} \ \qup\}.
 \end{equation}
Next, if  in a complex neighbourhood $U_c\subset\Con$ of $\qup$ there exists an analytic function $\zeta$ such that $\zeta(x)=\widehat{\mu}(x)$ for all $x\in\qup$, then $T_{\mu}$ coincides on $\bip$ with the differential operator $\zeta(D)$, where
\[
D=\Bigl(-i\frac{\partial}{\partial x_1},\dots, -i\frac{\partial}{\partial x_n}\bigr),
\]
(see \cite[p. 646]{3}). It is clearly that, given $m\in E^n$,  there exists $\mu_m\in \M$ such that
\[
\widehat{\mu}_m(t)=\prod_{j=1}^n t_j^{m_j}
\]
for all $t\in\qup$. Therefore, if $f\in\bip$, then
\begin{equation}
\partial^m_z f(x)= \widehat{\mu}_m(D)f(x)=\int_{\Rn}f(x-y)\,d\mu_m(y).
 \end{equation}
Assume  now that $f\in\bnn$. Fix any $y\in\Rn$. Then the function $f_y(x):=f(x-y)$, $x\in\Rn$,  is also in $\bnn$ Finally, using  that $\|\mu_m\|_{\M}<\infty$ and that $\bnn$ is a closed subspace of $\bip$,  we conclude from (2.6) that $\partial^m_z f\in\bnn$. The lemma is proved.

\ \ Recall that any $f\in B^{\infty}_{\qup}$ satisfies
\begin{equation}
|f(z)|\le \sup_{x\in\Rn}|f(x)|e^{\pi\sum_{j=1}^n |y_j|},
 \end{equation}
where $z=x+iy$, $x,y\in\Rn$ (see, e.g., \cite[p. 117]{10}). In addition, (2.1)  implies  that this estimate  is also true for any  $f\in \bnp$, $1\le p<\infty$ .
\begin{prop}
Let $f\in \bnn$. Suppose that
\begin{equation}
\partial^k_z f(u)= 0
 \end{equation}
for each $k\in E^n$ and all $u\in 2\Zn$. Then $f\equiv 0$.
\end{prop}
\ \ {\bf Proof}.\ Our proof is by induction on the dimension $n$ of $\Con$. If n = 1, then (2.8) is equivalent to
\begin{equation}
f(2k)=f'(2k)=0
 \end{equation}
for each $f\in B_{\pi}$ and all $k\in\Z$. Set
\begin{equation}
g(z)=\frac{f(z)}{\sin^2(\pi z/2)}.
 \end{equation}
Then (2.9) implies that $g$ is an entire function on $\Co$. For each $\tau\in\Z$, $\tau>0$, let us define
\[
D_{\tau}=\{z\in\Co: |\Re z|\le 1+2\tau, \ |\Im z|\le 2\}.
\]
Then
\begin{equation}
\min_{z\in \partial D_{\tau}}\Bigl|\sin\frac{\pi z}{2}\Bigr|\ge 1
 \end{equation}
for all $\tau=1,2,\dots$. Combining (2.7) with (2.11) and using the maximum modules principe for analytic function in the domain $D_{\tau}$, we see that $|g|$ is bounded on each $D_{\tau}$ by the same finite constant. Hence, $|g|$ is bounded on $D=:\cup_{p=1}^{\infty} D_{\tau}=\{z\in\Co: |\Im z|\le 2\}$. On the other hand, if $z\in\Co\setminus D$,  then it is easy to verify that
\begin{equation}
\Bigl|\sin\frac{\pi z}2\Bigr|\ge \frac12 e^{\pi|\Im z|/2|}.
 \end{equation}
Combining this with (2.7), we conclude that $|g|$ is bounded  on $\Co\setminus D$. Hence $ g$ is a constant $c\in\Co$. Then (2.10) gives  that $f(z)=c\sin^2(\pi z/2)$. Finally, $c=0$, since $\lim _{x\in\R; x\to\infty}f(x)=0$.

Suppose that  our proposition holds for dimension $m\ge 1$. First, we claim that if $f\in B_{Q^{m+1}_{\pi}}$  satisfies (2.8) with $n=m+1$, then
\begin{equation}
 \partial^s_z f(z)=0
 \end{equation}
for each $s\in E^{m+1}$ and all $z \in H^{m+1}$. Fix any $\widetilde{z}\in H^{m+1}$. According to (2.3), without loss of generality we can assume that $\widetilde{z}=(2u,\widetilde{z}_2,\dots,\widetilde{z}_n)$ with $u\in\Z$ and $(\widetilde{z}_2,\dots,\widetilde{z}_{m+1})\in \Co^m$. Let us define
\[
F(z)=f(2u, z_1,\dots,z_m)
\]
 for our fixed $u\in\Z$ and all $z=( z_1,\dots,z_m)\in\Co^m$. Then $F\in B_{Q^m_{\pi}}$. By the induction hypothesis,   the condition (2.8) with $n=m$   implies that
\[
 \partial^s_z F(z)=0
\]
for each $s\in E^{m}$ and all $z \in 2\Z^{m}$. Therefore,  the induction hypothesis gives $F\equiv 0$ on $\Co^m$, yielding our  claim (2.13).

\ \ Now from lemma 2.1 it follows that there is an entire function $g$ on $\Co^{m+1}$ such that
\begin{equation}
f(z)={\rm{sic}}_{m+1}(  \pi z ) g(z).
\end{equation}

\ \ Second, we claim that
\begin{equation}
g(z)=0
 \end{equation}
for all $z \in H^{m+1}$. Fix any $\widetilde{z} \in H^{m+1}$. We need to show that $g(\widetilde{z})=0$.  By (2.3), we can assume that there is   $r\in\Z$, $1\le  r\le m+1$, such that
\[
\widetilde{z}=(2u_1,\dots,2u_r,\widetilde{z}_{r+1},\dots,\widetilde{z}_{m+1})
\]
for  certain $u_1,\dots,u_r\in\Z$ and some $\widetilde{z}_{r+1},\dots,\widetilde{z}_{m+1}\not \in 2\Z$. Using our fixed numbers $u_1,\dots,u_r\in\Z$, let us define on $\Co^{m+1-r}$ the function
\begin{equation}
F_{\widetilde{z}}(\lambda)=\frac{\partial^{r}}{\partial z_1\cdots\partial z_r}f(2u_1,\dots,2u_r,\lambda_1,\dots,\lambda_{m+1-r}),
 \end{equation}
$\lambda\in \Co^{m+1-r}$. Then (2.8) implies that
\begin{equation}
\frac{\partial^{|\omega|}}{\partial^{\omega_1}\lambda_1\cdots\partial^{\omega_{m+1-r}}\lambda_{m+1-r}}F_{\widetilde{z}}(\lambda)=0
 \end{equation}
for each $\omega\in E^{m+1-r}$ and  all $\lambda\in 2\Z^{m+1-r}$. Next, Lemma  2.2 shows that $F_{\widetilde{z}}\in B_{Q^{m+1-r}_{\pi}}$. Therefore, using the induction hypothesis for dimension $m$ and keeping  in mind that $m+1-r\le m$, we conclude from (2.17) that $F_{\widetilde{z}}\equiv 0$ on  $\Co^{m+1-r}$. On the other hand,  (2.14) implies that
\[
F_{\widetilde{z}}(\lambda)=(-1)^{u_1+\dots+u_r}\Bigl(\frac{\pi}2\Bigr)^r\ \prod_{j=1}^{m+1-r} \sin\Bigl(\frac{\pi}2\lambda_j\Bigr)g(2u_,\dots,2u_r,\lambda_1,\dots,\lambda_{m+1-r}),
\]
for all $\lambda\in \Co^{m+1-r}$. Finally, if we take here $\lambda_j=\widetilde{z}_{r+j}$, $j=1,\dots,m+1-r$, use the fact that $F_{\widetilde{z}}\equiv 0$ on $\Co^{m+1-r}$, and keeping in mind that  $\widetilde{z}_{r+j}\not \in 2\Z$ for all $j=1,\dots,m+1-r$, then we get  $g(\widetilde{z})=0$. This proves (2.14).

\ \  Now lemma 2.1 shows that there is an entire function $h$ on $\Co^{m+1}$ such that $ g(z)={\rm{sic}}_{m+1}(\pi z) h(z)$. Using (2.13), we get
\[
   h(z)=\frac{f(z)}{{\rm{sic}}_{m+1}^2(\pi z)}
  \]
 for all $z\in\Co^{m+1}$. Now combining (2.7) with (2.11) and (2.12), we conclude that  $h $ is bounded on $\Co^{m+1}$. By Liouville's theorem, it follows that $h$ is a constant $c\in\Co$. Therefore, $ f(z)=c\cdot{\rm{sic}}_{m+1}^2(\pi z)$, $z\in\Co^{m+1}$. Using the fact that $f\in\bnn$, i.e., $\lim_{|x|\to\infty}f(x)=0$, we see that $c=0$, which completes the proof of Proposition 2.3.

\ \ For $1\le p<\infty$, let $l^p_n$ denote the usual Banach space of sequences of complex numbers $\{c_n\}_{u\in\Zn}$  such that $\sum_{u\in\Zn}|c_u|^p<\infty$. By Nikol'skii's inequality (\cite[p. 123]{10}),  for any $0<\theta<\infty$, there exists a finite constant $a=a(p, \sigma, \theta)$ such that
\[
\biggl(\sum_{u\in\Zn}|f(\theta u)|^p\biggr)^{1/p}\le a \|f\|_{\bnp}
\]
for all $f\in\bnp$. Combining this with Bernstein's inequality in $\bnp$ (see \cite[p. 116]{10}), we deduce that a similar estimate holds also for all derivatives $\partial ^k_z f$ with $k\in E^n$, i.e.,
\[
\biggl(\sum_{u\in\Zn}|\partial^k_z f(\theta u)|^p\biggr)^{1/p}\le A \|f\|_{\bnp}
\]
for each $0<\theta<\infty$, certain $A=A(k, p, \sigma, \theta)<\infty$ and all $f\in\bnp$. In particular, this  this means that if $f\in\bnp$, then
\begin{equation}
\Bigl\{\partial^k_z f(\theta u)\Bigr\}_{u\in \Zn}\in l^{p}_n
 \end{equation}
for each $0<\theta<\infty$ and all $k\in E^n$.

\ \ Recall that if $1<r<\infty$, then  (see \cite[p. 811]{11})
 \begin{equation}
 \sum_{m\in\Zn}\biggl|\sinc\Bigl(ax-m\Bigr)\biggr|^{r} \le\Bigl(1+\frac1{r-1}\Bigr)^n
\end{equation}
 for each $a>0$ and all $x\in\Rn$.

{\bf{Proof of Theorem 1.1.}} \ Let us  first show  that (1.7) converges absolutely and uniformly on $\Rn$. To this end, we divide (1.7) into $2^n$ series of the following type\begin{equation}
S_k(x)=\sum_{u\in 2\Zn}P_{f,k,u}(x-u)\sinc^2\Bigl(\frac{x-u}2\Bigr)=\sum_{u\in 2\Zn}\biggl(\partial^k_x f(u)\prod_{j=1}^n(x_j-u_j)^{k_j}\sinc^2\Bigl(\frac{x-u}2\Bigr)\biggr),
\end{equation}
$x\in\Rn$, $k\in E^n$.  Now it remains to prove that for any $\varepsilon>0$ there is a positive integer $\tau$ such that
\begin{equation}
|S_{k, \tau}(x)|=\Bigl|\sum_{\substack{u\in2\Zn\\ |u_1|,\dots|u_n|\ge \tau}}\partial^k_x f(u)\prod_{j=1}^n(x_j-u_j)^{k_j}\sinc^2\Bigl(\frac{x-u}2\Bigr)\Bigr|<\varepsilon
\end{equation}
for all $x\in\Rn$ and each $k\in E^n$.

 \ \ Fix an arbitrary number $p_1$ such that $p_1\ge p$ and $1< p_1<\infty$. Let $q_1=p_1/(p_1-1)$. Applying H\"{o}lder's inequality to (2.20)  gives
\begin{equation}
|S_{k, \tau}(x)|\le\Bigl(\sum_{\substack{u\in2\Zn\\ |u_1|,\dots|u_n|\ge \tau}}\Bigl|\partial^k_x f(u)\Bigr|^{p_1}\Bigr)^{1/p_1}\biggl(\sum_{\substack{u\in2\Zn\\ |u_1|,\dots|u_n|\ge \tau}}\prod_{j=1}^n|x_j-u_j|^{q_1k_j}\Bigl|\sinc\Bigl(\frac{x-u}2\Bigr)\Bigr|^{2q_1}\biggr)^{1/q_1}.
\end{equation}
 The condition $p_1\ge p$ implies that $l^p_n\subset l^{p_1}_n$. Hence, we conclude from (2.18) that there exists a positive integer $\tau_1$ such that
\begin{equation}
\Bigl(\sum_{\substack{u\in2\Zn\\ |u_1|,\dots|u_n|\ge \tau_1}}\Bigl|\partial^k_x f(u)\Bigr|^{p_1}\Bigr)^{1/p_1}<\varepsilon.
\end{equation}
Since $|\sivc t|\le 1$ for $t\in\R$ and $k_1,\dots,k_n\in\{0;1\}$, we see that
\begin{gather}
|x_j-u_j|^{q_1k_j}\Bigl|\sivc\Bigl(\frac{x_j-u_j}2\Bigr)\Bigr|^{2q_1}\le \Bigl(\frac2{\pi}\Bigr)^{q_1k_j}\Bigl|\sivc\Bigl(\frac{x_j-u_j}2\Bigr)\Bigr|^{q_1(2-k_j)}\nonumber \\
\le \Bigl(\frac2{\pi}\Bigr)^{q_1k_j}\Bigl|\sivc\Bigl(\frac{x_j-u_j}2\Bigr)\Bigr|^{q_1}\nonumber
\end{gather}
for all $x_j\in\R$ and  each $u_j\in\Z$, $ j=1,\dots,n$. Therefore, the second factor on the right of (2.22) is  estimated by
\begin{gather}
\sum_{\substack{u\in2\Zn\\ |u_1|,\dots|u_n|\ge \tau}}\prod_{j=1}^n|x_j-u_j|^{q_1k_j}\Bigl|\sinc\Bigl(\frac{x-u}2\Bigr)\Bigr|^{2q_1}\nonumber \\
\le \Bigl(\frac{2}{\pi}\Bigr)^{(k_1+\dots+k_n)q_1}
\sum_{\substack{u\in2\Zn\\ |u_1|,\dots|u_n|\ge \tau}}\prod_{j=1}^n\Bigl|\sivc\Bigl(\frac{x_j-u_j}2\Bigr)\Bigr|^{q_1}\nonumber \\
= \Bigl(\frac{2}{\pi}\Bigr)^{(k_1+\dots+k_n)q_1}\sum_{\substack{u\in2\Zn\\ |u_1|,\dots|u_n|\ge \tau}}\Bigl|\sinc\Bigl(\frac{x-u}2\Bigr)\Bigr|^{q_1},\nonumber \\
\end{gather}
$x\in\Rn$. Now the condition $1<p_1<\infty$ implies that $1<q_1<\infty$. Therefore, using (2.19), we conclude that there exists positive integer $\tau_2$ such that
\begin{equation}
\biggl(\sum_{\substack{u\in2\Zn\\ |u_1|,\dots|u_n|\ge \tau_2}}\prod_{j=1}^n|x_j-u_j|^{q_1k_j}\Bigl|\sinc\Bigl(\frac{x-u}2\Bigr)\Bigr|^{2q_1}\biggr)^{1/q_1}<\varepsilon
\end{equation}
for all $x\in\Rn$. Finally, if we assume that $\varepsilon<1$ and take $\tau=\max\{\tau_1, \tau_2\}$, then  combining (2.22), (2.23) and (2.25), we obtain (2.21). This complete the proof that (1.7)  converges absolutely and uniformly on $\Rn$.

\ \ If $f\in\bnp$, $1< p<\infty$, then any partial sum of (1.7)  is also  in $\bnp$. For $f\in B^1_{Q^n_{\pi}}$, these partial sums  is  in $\bnp$ for each $1< p<\infty$. Hence, by (2.1),  these sums are also  elements of $\bnn$ for any $1\le p<\infty$.   Let $F$ be the sum of  (1.7). Since $\bnn$ is a Banach space and  (1.7)  converges absolutely and uniformly on $\Rn$, it follows that $F\in\bnn$. Let us denote $f_1=f-F$. Then $f_1\in\bnn$. Using the fact that (1.7) converges  uniformly on $\Rn$ to $F$, we conclude that  $f_1$ satisfies all conditions (2.8). By Proposition 2.3, we have that $f_1\equiv 0$. Hence, $F\equiv f$, i.e., equality (1.7)  holds for all $x\in\Rn$.

\ \ Let $K$ be a compact subset of $\Con$. Then there exists $0<a<\infty$ such that $K$ is a subset of the strip $M_a=\{z\in\Con: |\Im z_j|\le a, j=1,\dots,n\}$. Given $\omega\in\Zn$, let us define the partial sum of (1.7) by
\[
W_{\omega}(z)=\sum_{\substack{u\in2\Zn\\ |u_1|\ge|\omega_1|,\dots|u_n|\ge |\omega_n|}}\biggl(\sum_{k\in E^n}P_{f,k,u}(z-u)\biggr)\sinc^2\Bigl(\frac1{2\pi}\sigma(z-u)\Bigr).
\]
According to (2.1), we have from (2.7) that
\begin{equation}
\sup_{z\in K}|f(z)-W_{\omega}(z)|\le \sup_{x\in \Rn}|f(x)-W_{\omega}(x)|e^{\pi a n}.
\end{equation}
Since (1.7)  converges uniformly on $\Rn$ to $f$, we see from (2.26) that (1.7) converges uniformly also  on $K$ to $f$. Theorem 1.1 is proved.

}}
\end{document}